\documentclass[10pt]{article}
\newcommand{\D}{\displaystyle}
\newcommand{\DF}[2]{\frac{\D#1}{\D#2}}
\usepackage{CJK}
\usepackage{amsmath}
\usepackage{amssymb}
\usepackage{listings}
\usepackage{newlfont}
\usepackage{algorithm}
\usepackage{multirow}
\usepackage{graphicx}
\usepackage{subeqnarray}
\usepackage{cases}
\usepackage{latexsym,bm}
\usepackage{indentfirst}
\usepackage{hyperref}
\usepackage{epstopdf}

\usepackage{color}

\numberwithin{equation}{section}

\date{}
\title{A primal-dual fixed point algorithm for multi-block convex minimization}
\author{Peijun Chen$^{1,2,3}$,\  Jianguo Huang$^{1}$,\ Xiaoqun Zhang$^{1,4}$
\\
\\$^{1}${\small \textit{School of Mathematical Sciences, and MOE-LSC, Shanghai Jiao Tong University,}} \\{\small \textit{Shanghai 200240, China}}
\\
$^{2}${\small \textit{School of Biomedical Engineering, Shanghai Jiao Tong University,}}\\{\small \textit{Shanghai 200240, China}}
\\
$^{3}${\small \textit{ Department of Mathematics, Taiyuan University of Science and Technology,}}
\\{\small \textit{Taiyuan 030024, China}}
\\$^{4}$ {\small\textit{Institute of Natural Sciences, Shanghai Jiao Tong University, }}
\\
{\small\textit{Shanghai 200240, China}}
\\
{\small Email: chenpeijun@sjtu.edu.cn, jghuang@sjtu.edu.cn and xqzhang@sjtu.edu.cn}
}

\begin{document}
\maketitle

\noindent
\begin{abstract}
We extend a primal-dual fixed point algorithm (PDFP) proposed in \cite{CHZ15}
to solve two kinds of separable multi-block minimization problems, arising in signal processing and imaging science. This work shows the flexibility of applying PDFP algorithm to multi-block problems and illustrate how practical and fully decoupled schemes can be derived, especially for  parallel implementation of  large scale problems. The connections and comparisons to  the alternating direction method of multiplier (ADMM) are also present. We demonstrate how different algorithms can be obtained by splitting the problems in different ways through the classic example of  sparsity regularized least square  model with constraint.  In particular, for a class of linearly constrained problems, which are of great interest in the context of multi-block ADMM, can be solved by PDFP with a guarantee of convergence.  Finally, some experiments are provided to illustrate the performance of several schemes derived by the PDFP algorithm.
\end{abstract}

\section{Introduction}
In this paper, we are concerned with extending the primal-dual fixed point (PDFP) algorithm proposed in \cite{CHZ15} for solving two kinds of general multi-block problems \eqref{eqbasic_TVL1} and \eqref{eqbasic_MSB} with maximally decoupled iterative scheme. The first kind of problems are formulated as
\begin{equation}
   \underset{x\in  \mathbb{R}^n}{\mbox{min}}\quad  {f_1}(x)+\sum_{i=1}^N\theta_i(B_ix+b_i)+f_3(x), \label{eqbasic_TVL1}
\end{equation}
where $\theta_i\in \Gamma_0 (\mathbb{R}^{m_i})$,  $B_i:{\mathbb{R}^n}\rightarrow \mathbb{R}^{m_i}$ a  bounded linear transform, $b_i\in \mathbb{R}^{m_i}$, $i=1,2,\cdots,N$.  ${{f_1}},{f_3}\in \Gamma_0 (\mathbb{R}^n)$  and ${{f_1}}$ is differentiable on  ${\mathbb{R}^n}$ with  $1/\beta$-Lipschitz continuous gradient for some $\beta\in(0,+\infty]$.  Throughout this paper, $\Gamma_0(\mathbb{R}^n)$ stands for  the collection of all proper lower semicontinuous convex functions from $\mathbb{R}^n$ to $(-\infty,+\infty]$. Many problems in image processing and signal recovery with multi-regularization terms can be formulated in the form of \eqref{eqbasic_TVL1}.

The second kind of problems under discussion are optimization problems with constraints, given as follows.
\begin{equation}\label{eqbasic_MSB}
\begin{aligned}
   &\underset{x_1,x_2,\cdots,x_N}{\mbox{\quad min \quad}} \sum_{i=1}^{N_1}\theta_i(B_ix_i+b_i)+\sum_{i=N_1+1}^N\theta_i(x_i)\\
   &\mbox{\quad\  st. }\sum_{i=1}^N A_ix_i=a,\\
   &\phantom{\mbox{\quad\  st. }} x_i\in C_i,i=1,2,\cdots,N.
\end{aligned}
\end{equation}
Here,  $\theta_i\in \Gamma_0 (\mathbb{R}^{m_i})$, $B_i:{\mathbb{R}^{n_i}}\rightarrow \mathbb{R}^{m_i}$   a bounded linear transform and $b_i\in \mathbb{R}^{m_i}$ for  $i=1,2,\cdots,N_1$. Moreover, for $i=N_1+1,\cdots,N$, $\theta_i\in \Gamma_0 (\mathbb{R}^{n_i})$ is differentiable on  $\mathbb{R}^{n_i}$ with  $1/\beta_i$-Lipschitz continuous gradient for some $\beta_i\in(0,+\infty]$. For $i=1,2,\cdots,N$, the  constraint set $C_i\subset \mathbb{R}^{n_i}$ is closed and convex, $A_i$ is a $l\times n_i$ matrix, and $a\in \mathbb{R}^{l}$.

Many problems can be formulated in the form \eqref{eqbasic_MSB}, for example elliptic optimal control problems \cite{CK11}. In some applications, the problem \eqref{eqbasic_TVL1} can be viewed as a decomposition on the observed data, while the problem \eqref{eqbasic_MSB} is a mixture of the variables and data decomposition. In particular, for some special cases, both problems
\eqref{eqbasic_TVL1} and \eqref{eqbasic_MSB} can be abstracted as
\begin{equation}\label{eqbasic_MSB02}
\begin{aligned}
   &\underset{x_1,x_2,\cdots,x_N}{\mbox{\quad min \quad}} \sum_{i=1}^{N}\theta_i(x_i)\\
   &\mbox{\quad\  st. }\sum_{i=1}^N A_ix_i=a,\\
   &\phantom{\mbox{\quad\  st. }} x_i\in C_i,i=1,2,\cdots,N,
\end{aligned}
\end{equation}
by properly introducing auxiliary variables, or vice-visa, depending on
the simplicity of the functions $\theta_i$ involved. In the literature, many existing works have been devoted to solving \eqref{eqbasic_MSB02}, for example, the variants of popular alternating direction
method of multipliers (ADMM) \cite{HTY12,HXY13,DLP13} for three or more block problems.

Now, let us recall the proximal primal-dual fixed point algorithm {PDFP} in \cite{CHZ15}
for solving the following three-block problem
\begin{equation}
\underset{x\in \mathbb{R}^n}{\mbox{min}}\quad {{f_1}}(x)+ {{f_2}}(Bx+b)+f_3(x). \label{PDFP2O3B:eqbasic}
\end{equation}
In \eqref{PDFP2O3B:eqbasic}, ${{f_2}}\in \Gamma_0 (\mathbb{R}^m)$, $B:{\mathbb{R}^n}\rightarrow \mathbb{R}^m$ a bounded linear transform, $b\in \mathbb{R}^m$, ${{f_1}}$ and ${f_3}$ are the same ones as given in \eqref{eqbasic_TVL1}. As usual, define the proximity operator $\mbox{prox}_f$ of $f$ by (cf. \cite{CW05})
\[
   \mbox{prox}_{f}(x)= \underset{y\in\mathbb{R}^n}{\mbox{arg min}}\ {f
    (y)+\DF{1}{2}\|x-y\|^2}.
\]
Then, our PDFP algorithm can be described as follows.
\begin{subequations}
 \label{formbasic3B}
 \begin{numcases} {(\mbox{PDFP} ) \quad }
   x^{k+1/2}=\mbox{prox}_{{\gamma}{f_3}}(x^k-\gamma\nabla {{f_1}}(x^k)-{\lambda} B^T  v^{k}),\label{formbasic3Ba}\\
   v^{k+1}=(I-\mbox{prox}_{\frac{\gamma}{\lambda}{{f_2}}})(Bx^{k+1/2}+b+v^k)\label{formbasic3Bb},\\
   x^{k+1}=\mbox{prox}_{{\gamma}{f_3}}(x^k-\gamma\nabla {{f_1}}(x^k)-{\lambda} B^T  v^{k+1}),\label{formbasic3Bc}
 \end{numcases}
\end{subequations}
where $0<\lambda< 1/\lambda_{\max}(BB^T)$ and $0<\gamma<2\beta$.

The purpose of this paper is intended to extend {PDFP} to solve the above two kinds of general multi-block problems \eqref{eqbasic_TVL1} and \eqref{eqbasic_MSB} with maximally decoupled iteration scheme. The key trick of our treatment is the use of {PDFP} combined with feasible reformulation of the multi-block problems in the form \eqref{PDFP2O3B:eqbasic}, so that we can derive many variants of iterative schemes with different structures. One obvious advantage of the extended schemes is their simplicity and the convenience for parallel implementation. Some of the algorithms derived in this paper already exist in the literature and some of them are new and effective. The new schemes are compared with the ADMM and we will show the connection and the difference later on. We mention in passing that similar techniques are also adopted in \cite{C13,CHZ1302,LSXZ15,TZW15}. Compared to the schemes developed in \cite{C13,LSXZ15,TZW15}, if a scheme is established based on PDFP with $f_1$ nonzero in \eqref{PDFP2O3B:eqbasic}, it's more convenient for us to choose parameters in applications, as shown in \cite{CHZ15}. However, if a scheme is constructed based on PDFP  by viewing $f_1$ equal to $0$, it requires to compute an additional symmetric step. Note that in many $\ell_1$-based regularization problems, this step can be implemented explicitly. So the additional cost is ignorable. In what follows, to simplify the presentation, we will not systematically compare the schemes developed here with those in \cite{CP11,C13,LSXZ15,LZ15} any more.

The rest of the paper is organized as follows. In Section \ref{sec:mblock}, we will show how PDFP can be extended to solve \eqref{eqbasic_TVL1}, present the connections and differences with ADMM and derive different algorithms by using the constrained and sparse regularized image restoration  model as an illustrative  example. In Section \ref{sec:MSB},  {PDFP} is extended to solve \eqref{eqbasic_MSB}, and we also show the comparison with ADMM.  In Section \ref{sec:numericalE}, the numerical performance and efficiency of the variants of {PDFP}  are demonstrated through constrained total variation computerized tomography (CT) reconstruction and solving quadratic programming model.

\section{PDFP for the muti-block problem \eqref{eqbasic_TVL1}}\label{sec:mblock}
\subsection{Algorithm and its deduction}\label{sec:alg_mblock}
In this section, we formulate \eqref{eqbasic_TVL1} as  a special case of \eqref{PDFP2O3B:eqbasic}. Then the PDFP algorithm can be applied and formulated in parallel form due to  the separability of ${f_2}$ on its variables. Similar technique has also been used in \cite{C13,CHZ1302,LSXZ15,TZW15} and we present the details here for completeness.

Rewrite the second term in \eqref{eqbasic_TVL1}  as
\begin{align*}
  {f_2}(Bx+b):=\sum_{i=1}^{N}\theta_i(B_ix+b_i)
\end{align*}
with the symbols
\begin{gather*}
f_2(y)=\sum_{i=1}^{N}\theta_i(y_i), y=Bx+b,
\\
    B=
  \begin{pmatrix}
       B_1\\
       B_2\\
       \vdots\\
       B_N
       \end{pmatrix},
\
  b=
  \begin{pmatrix}
       b_1\\
       b_2\\
       \vdots\\
       b_N
       \end{pmatrix}.
  \end{gather*}
Thus, the problem \eqref{eqbasic_TVL1} can be recast in the form of \eqref{PDFP2O3B:eqbasic} and is resolved with PDFP. Since
${f_2}$ is separable in terms of its variables, the scheme \eqref{formbasic3B} can be further expressed as
\begin{subequations}
 \label{formbasicV}
 \begin{numcases} { }
   x^{k+1/2}=\mbox{prox}_{{\gamma}{f_3}}(x^k-\gamma \nabla {{f_1}}(x^k)-\lambda \sum_{j=1}^{N}B_j^Tv_j^{k}),\label{formbasicVa}\\
   v_i^{k+1}=(I-\mbox{prox}_{\frac{\gamma}{\lambda}{\theta_i}})(B_ix^{k+1/2}+b_i+v_i^k),i=1,2,\cdots,N\label{formbasicVb},\\
   x^{k+1}=\mbox{prox}_{{\gamma}{f_3}}(x^k-\gamma \nabla {{f_1}}(x^k)-\lambda \sum_{j=1}^{N}B_j^Tv_j^{k+1}).\label{formbasicVc}
 \end{numcases}
\end{subequations}
The convergence condition of PDFP in \cite{CHZ15} implies that the above algorithm is convergent whenever  $0<\lambda< 1/\sum_{i=1}^{N}\lambda_{\max}(B_iB_i^T)$ and $0<\gamma<2\beta$. The scheme \eqref{formbasicV} is naturally in a parallel form,  which may be useful for large scale problems. Also for some special cases, such as $f_1=0$, $f_3=\chi_C$,
one may even get simpler forms (see \cite{CHZ15} for details).

\subsection{Comparison to ADMM}\label{sec:CmpWithADMM}

There are  many works on ADMM methods \cite{HTY12,HXY13,DLP13}. We will show the difference between PDFP and ADMM for solving \eqref{eqbasic_TVL1}. Since our method for solving \eqref{eqbasic_TVL1} is based on the PDFP \eqref{formbasic3B} for solving \eqref{PDFP2O3B:eqbasic}. We first show how the ADMM resolves the same problem. In fact, we should first reformulate the problem in the form \eqref{eqbasic_MSB02} by introducing auxiliary variables. Then, we can use the ADMM to drive the scheme for solving \eqref{PDFP2O3B:eqbasic}. However, our PDFP is developed based on a fixed point formulation the solution of \eqref{PDFP2O3B:eqbasic} must satisfy. So the ideas of constructing the two methods are quite different.

To show the difference of the two methods more clearly, we compare their schemes for solving \eqref{eqbasic_TVL1} with $f_3=0$. PDFP for solving \eqref{eqbasic_TVL1} have been given in \eqref{formbasicV} based on three blocks algorithm \eqref{formbasic3B}. We can also use the similar technique to achieve the ADMM method in this case:
\begin{subequations}
 \label{_MV_ADMM}
 \begin{numcases} { }
   x^{k+1}=\underset{x\in \mathbb{R}^n}{\mbox{argmin  }} f_1(x)+\frac {\beta} 2\sum_{i=1}^N\|B_ix+(b_i-y_i^k+v_i^k)\|^2,\label{_MV_ADMMa}\\
   y_i^{k+1}=\mbox{prox}_{\frac{1}{\beta}{\theta_i}}(B_ix^{k+1}+b_i+v_i^k),i=1,2,\cdots,N,\label{_MV_ADMMb}\\
   v_i^{k+1}=v_i^k+ \tau (B_ix^{k+1}+b_i-y_i^{k+1}),i=1,2,\cdots,N. \label{_MV_ADMMc}
  \end{numcases}
\end{subequations}
As a matter of fact, the scheme \eqref{_MV_ADMM} follows from an application of the two block ADMM for solving \eqref{eqbasic_TVL1}  with $f_3=0$ and the convergence condition for \eqref{_MV_ADMM} is still $\beta>0$  and $\tau\in (0,(1+\sqrt{5})/2)$.

\subsection{Application to constrained sparse regularization problems}
In this subsection, we will consider how to get different algorithms by using the extension of PDFP \eqref{formbasicV} for a specific problem. The problem that we are interested is the well-known constrained  sparse regularization model in inverse problems and imaging:
\begin{equation}
   \underset{x\in {C}}{\mbox{ min}}\quad \frac 1 2\|Ax-a\|^2+\mu \|Dx\|_1, \label{eq_TVL2}
\end{equation}
where $\|Ax-a\|^2$ is the smooth data-fidelity term, $\mu \|Dx\|_1$ is the regularization term to ensure the solution is sparse under the transform $D$ and $\mu$ is the regularization parameter. The problem \eqref{eq_TVL2} is equivalent to
\begin{equation}
   \underset{x\in \mathbb{R}^n}{\mbox{min}}\quad \frac 1 2\|Ax-a\|^2+\mu \|Dx\|_1+\chi_C(x), \label{eq_TVL202}
\end{equation}
where \begin{equation*}
   \chi_C(x)=
     \left\{
      \begin{array}{ll}
         0, &x\in C,\\
         +\infty,&x\not\in C.
      \end{array}
     \right.
\end{equation*}

First, applying PDFP  \eqref{formbasic3B} to the problem \eqref{PDFP2O3B:eqbasic} with the three blocks given by $f_1(x)=\frac 1 2\|Ax-a\|^2$, $f_2=\mu\|\cdot\|_1$, $B=D$, $b=0$, $f_3=\chi_C$ and  noting $\mbox{prox}_{{\gamma}{\chi_C}}=\mbox{proj}_C$, we obtain
\begin{align}
 \label{formTVL2C2}
 \mbox{(Scheme 1) } \left\{
 \begin{aligned}
   &x^{k+1/2}=\mbox{proj}_C(x^k-\gamma A^T(A x^k-a)-{\lambda} D^T  v^{k}),\\
   &v^{k+1}=(I-\mbox{prox}_{\frac{\gamma}{\lambda}{\mu\|\cdot\|_1}})(Dx^{k+1/2}+v^k),  \\
   &x^{k+1}=\mbox{proj}_C(x^k-\gamma A^T(A x^k-a)-{\lambda} D^T  v^{k+1}),
 \end{aligned}
 \right.
\end{align}
where $0<\lambda< 1/ \lambda_{\max}(DD^T) $ and $0<\gamma<2/ \lambda_{\max}(A^TA)$. This is the original algorithm proposed in \cite{CHZ15}.

The second scheme can be obtained  by setting $f_1(x)=\frac 1 2\|Ax-a\|^2$, $\theta_1=\mu\|\cdot\|_1$, $B_1=D$, $b_1=0$, $\theta_2=\chi_C$, $B_2=I$, $b_2=0$, $f_3=0$, leading to
\begin{align}
 \label{formTVL2C1}
 \mbox{(Scheme 2) }\left\{
 \begin{aligned}
   &x^{k+1/2}=x^k-\gamma  A^T(A x^k-a)-{\lambda} D^T  v_1^{k}-{\lambda}v_2^{k}, \\
   &v_1^{k+1}=(I-\mbox{prox}_{\frac{\gamma}{\lambda}{\mu\|\cdot\|_1}})(Dx^{k+1/2}+v_1^k),  \\
   &v_2^{k+1}=(I-\mbox{proj}_{C})(x^{k+1/2}+v_2^k), \\
   &x^{k+1}=x^k-\gamma A^T(A x^k-a)-{\lambda} D^T  v_1^{k+1}-{\lambda}v_2^{k+1},
 \end{aligned}
 \right.
\end{align}
where $0<\lambda\leq 1/ (\lambda_{\max}(DD^T)+1) $ and $0<\gamma<2/ \lambda_{\max}(A^TA)$. This scheme \eqref{formTVL2C1} is the form proposed in \cite{CHZ1302} by recasting the problem in two-block. We note that $x^{k+1}$ may not be a feasible solution during the iteration. In addition, an auxiliary variable $v_2$ is introduced and the permitted ranges of the parameter $\lambda$ is also a little tighter compared to Scheme 1.

In the following,  we present some schemes to use different properties of the objective functions $\frac 1 2\|Ax-a\|^2$, which may be the main computation cost in inverse problem applications.  By setting $f_1(x)=0$, $\theta_1=\mu\|\cdot\|_1$, $B_1=D$, $b_1=0$, $\theta_2=\frac 1 2\|Ax-a\|^2$, $B_2=I$, $b_2=0$, $f_3=\chi_C$, we can use
 \eqref{formbasicV} to solve \eqref{eq_TVL2} and obtain
\begin{align}
 \label{formTVL2nof1}
 \mbox{}\left\{
 \begin{aligned}
   &x^{k+1/2}=\mbox{proj}_C(x^k- \lambda(D^Tv_1^k+v_2^k)),\\
   &v_1^{k+1}=(I-\mbox{prox}_{\frac{\gamma}{\lambda}{\mu\|\cdot\|_1}})(Dx^{k+1/2}+v_1^k),\\
   &v_2^{k+1}=(x^{k+1/2}+v_2^k)-(I+{\frac{\gamma}{\lambda}A^TA})^{-1}(\frac{\gamma}{\lambda}A^Ta+x^{k+1/2}+v_2^k),\\
   &x^{k+1}=\mbox{proj}_C(x^k- \lambda(D^Tv_1^{k+1}+v_2^{k+1})),
 \end{aligned}
 \right.
\end{align}
where $0<\lambda< 1/ (\lambda_{\max}(DD^T)+1)$ and $0<\gamma<+\infty $. This scheme can be practical when the inverse of the matrix $(I+{\frac{\gamma}{\lambda}A^TA})$ is easy to obtain, for examples, for some diagonalizable matrix $A^TA$.

When the inverse of the matrix  $(I+{\frac{\gamma}{\lambda}A^TA})$ is not easy to compute, we can rewrite $\frac 1 2 \|Ax-a\|^2$ as $\frac 1 2 \|
\cdot\|^2\circ (Ax-a)$ and  set $f_1(x)=0$, $\theta_1=\mu\|\cdot\|_1$, $B_1=D$, $b_1=0$, $\theta_2=\frac{1}{2}\|\cdot\|^2$, $B_2=A$, $b_2=-a$, $f_3=\chi_C$,  and obtain
\begin{equation}
 \label{formTVL2nof102}
 \mbox{(Scheme 3) }\left\{
 \begin{aligned}
   &x^{k+1/2}=\mbox{proj}_C(x^k- \lambda(D^Tv_1^k+A^Tv_2^k)),\\
   &v_1^{k+1}=(I-\mbox{prox}_{\frac{\gamma}{\lambda}{\mu\|\cdot\|_1}})(Dx^{k+1/2}+v_1^k),\\
   &v_2^{k+1}=\frac{\gamma}{\gamma +\lambda}(Ax^{k+1/2}-a+v_2^k),\\
   &x^{k+1}=\mbox{proj}_C(x^k- \lambda(D^Tv_1^{k+1}+A^Tv_2^{k+1})),
 \end{aligned}
 \right.
\end{equation}
where $0<\lambda< 1/ (\lambda_{\max}(DD^T)+\lambda_{\max}(A^TA))$ and $0<\gamma<+\infty $.

If  we partition $A$ and $a$ into $N$ block rows, namely $A=(A_1^T,A_2^T,\cdots,A_N^T)^T$, $a=(a_1^T,a_2^T,\cdots,a_N^T)^T$,  where $A_j$ is a $m_j\times n$ matrix and $m_j<<n$ , $a_j\in \mathbb{R}^{m_j}$,  then $\DF{1}{2} \|A x -a\|^2=\DF{1}{2}\sum_{j=1}^N\|A_j x -a_j\|^2$. Here $A_i$ is different from the ones in \eqref{eqbasic_MSB}-\eqref{eqbasic_MSB02}, and they are  only used in this subsection.  It is very easy to see that  the scheme \eqref{formTVL2nof102} can be  written  in  a parallel form as
\begin{align}
 \label{formTVL2nof104}
 \left\{
 \begin{aligned}
   &x^{k+1/2}=\mbox{proj}_C(x^k- \lambda(D^Tv_1^k+\sum_{j=1}^NA_j^Tv_{2j}^k)),\\
   &v_1^{k+1}=(I-\mbox{prox}_{\frac{\gamma}{\lambda}{\mu\|\cdot\|_1}})(Dx^{k+1/2}+v_1^k),\\
   &v_{2i}^{k+1}=\frac{\gamma}{\gamma +\lambda}(A_ix^{k+1/2}-a_i+v_{2i}^k),i=1,2,\cdots,N,\\
   &x^{k+1}=\mbox{proj}_C(x^k- \lambda(D^Tv_1^{k+1}+\sum_{j=1}^NA_j^Tv_{2j}^{k+1})),
 \end{aligned}
 \right.
\end{align}
where $0<\lambda< 1/ (\lambda_{\max}(DD^T)+\sum_{j=1}^N \lambda_{\max}(A_j^TA_j))$ and $0<\gamma<+\infty $.

The above schemes except \eqref{formTVL2nof1} are fully explicit and involves only matrix-vector multiplication. In the following, we derive a semi-implicit scheme, which only involves the inverse of small size matrix. By setting $f_1(x)=0$, $\theta_1(x)=\mu\|x\|_1$, $B_1=D$, $b_1=0$, $\theta_{i+1}(x)=\DF{1}{2}\|A_ix-a_i\|^2$, $B_{i+1}=I$, $b_{i+1}=0$, $i=1,2,\cdots,N$, $f_3=\chi_C$, we obtain the following scheme by applying \eqref{formbasicV}:
\begin{equation}
 \label{formTVL2nof103}
 \mbox{(Scheme 4) }\left\{
 \begin{aligned}
   &x^{k+1/2}=\mbox{proj}_C(x^k- \lambda(D^Tv_1^k+\sum_{j=1}^N v_{2j}^{k+1})),\\
   &v_1^{k+1}=(I-\mbox{prox}_{\frac{\gamma}{\lambda}{\mu\|\cdot\|_1}})(Dx^{k+1/2}+v_1^k),\\
   &v_{2i}^{k+1}=(x^{k+1/2}+v_{2i}^k)-(I+{\frac{\gamma}{\lambda}A_i^TA_i})^{-1}(\frac{\gamma}{\lambda}A_i^Ta_i+x^{k+1/2}+v_{2i}^k), i=1,2,\cdots,N,\\
   &x^{k+1}=\mbox{proj}_C(x^k- \lambda(D^Tv_1^{k+1}+\sum_{j=1}^N v_{2j}^{k+1})),
 \end{aligned}
 \right.
\end{equation}
where $0<\lambda< 1/ (\lambda_{\max}(DD^T)+N)$ and $0<\gamma<+\infty $. At first glance, the size of the inverse in the third equation in \eqref{formTVL2nof103} is the same with the third ones in \eqref{formTVL2nof1}. However, thanks to the well known Sherman-Morrison-Woodbury formula, we know
\begin{align}
 \label{keytechnique}
 (I+{\frac{\gamma}{\lambda}A_i^TA_i})^{-1}=I-\frac\gamma\lambda A_i^T(I+\frac\gamma\lambda  A_iA_i^T)^{-1}A_i.
\end{align}
so we only need to invert a smaller size matrix  $I+\frac\gamma\lambda A_iA_i^T$ instead of $I+\frac\gamma\lambda A_i^TA_i$.
By using \eqref{keytechnique}, the third equation in \eqref{formTVL2nof103} is equivalent to
\begin{align}
\label{keytechnique02}
   v_{2i}^{k+1}=\frac\gamma\lambda A_i^T(I+\frac\gamma\lambda  A_iA_i^T)^{-1}A_i(\frac{\gamma}{\lambda}A_i^Ta_i+x^{k+1/2}+v_{2i}^k)-\frac{\gamma}{\lambda}A_i^Ta_i, i=1,2,\cdots,N.
\end{align}

\section{{PDFP} for constrained muti-block problem \eqref{eqbasic_MSB}}\label{sec:MSB}

In this section, we will show how to extend {PDFP} to solve \eqref{eqbasic_MSB}.  \eqref{eqbasic_MSB} can be also seen as a special case of \eqref{PDFP2O3B:eqbasic} by using operator $B$ and vector $b$, so we can solve it with {PDFP}. As a matter of fact, by using the separability of ${f_2}$ and ${f_1}$ about their variants, respectively, and noting that $C$ is separable, we can get the primal-dual fixed point algorithm \eqref{formbasicMC} for solving \eqref{eqbasic_MSB}.

\subsection{Algorithms and its deduction}\label{subsec:msb}
As a special case of indicator function $\chi_C$ on convex set $C$,
for $C=\{0\}$,  we define
\begin{equation*}
   \chi_0(x)=
     \left\{
      \begin{array}{ll}
         0, &x=0,\\
         +\infty,&x\neq 0.
      \end{array}
     \right.
\end{equation*}
Then \eqref{eqbasic_MSB} is equivalent to
\begin{equation}
   \underset{x_1,x_2,\cdots,x_N}{\mbox{  min }} \sum_{i=N_1+1}^N\theta_i(x_i)+\left(\sum_{i=1}^{N_1}\theta_i(B_ix_i+b_i)+\chi_{0}(\sum_{i=1}^N A_ix_i-d)\right)+\sum_{i=1}^N\chi_{C_i}(x_i) . \label{eqbasic_MSB2}
\end{equation}
Let
\begin{align*}
  &{f_1}(x)={f_1}(x_1,x_2,\cdots,x_N)=\sum_{i=N_1+1}^N\theta_i(x_i), \\
  &{f_3}(x)={f_3}(x_1,x_2,\cdots,x_N)=\sum_{i=1}^N\chi_{C_i}(x_i).
\end{align*}
Let
\begin{align*}
  &y_i=B_ix_i+b_i,i=1,2,\cdots,N_1,\\
  &y_{N_1+1}=\sum_{i=1}^NA_ix_i-a,
\end{align*}
\begin{align*}
  y=
  \begin{pmatrix}
       y_1\\
       y_2\\
       \vdots\\
       y_{N_1}\\
       y_{N_1+1}
       \end{pmatrix},
  B=
  \begin{pmatrix}
       B_1\\
       &B_2\\
       &&\ddots\\
       &&&B_{N_1}\\
       A_1& A_2&\cdots& A_{N_1}&\cdots&\ A_N
  \end{pmatrix},
  x=
  \begin{pmatrix}
       x_1\\
       x_2\\
       {  \vdots}\\
       x_{N_1}\\
       { \vdots}\\
       x_N
       \end{pmatrix},
  b=
  \begin{pmatrix}
       b_1\\
       b_2\\
       \vdots\\
       b_{N_1}\\
       -a
  \end{pmatrix},
\end{align*}
\begin{align*}
   {f_2}(y)={f_2}(y_1,y_2,\cdots,y_{N_1},y_{N_1+1})=\sum_{i=1}^{N_1}\theta_i(y_i)+\chi_{0}(y_{N_1+1}).
\end{align*}
Then we have
\begin{align*}
  y=Bx+b,\mbox{  } {f_2}(Bx+b)=\sum_{i=1}^{N_1}\theta_i(B_ix_i+b_i)+\chi_{0}(\sum_{i=1}^NA_ix_i-a),\label{notation_f1B}
\end{align*}
and problem \eqref{eqbasic_MSB} can be viewed as a special case of problem \eqref{PDFP2O3B:eqbasic}. Hence, we can use {PDFP} for solving \eqref{eqbasic_MSB}.
Observing that ${f_2}$ is separable about its variables $y_1,y_2,\cdots,y_{N_1+1}$,
${f_1}$ and $f_3$ are separable about their variables $x_1,x_2,\cdots,x_{N}$,  $\mbox{prox}_{\frac{\gamma}{\lambda}\chi_{C_i}}=\mbox{proj}_{C_i}$,
$\mbox{prox}_{\frac{\gamma}{\lambda}\chi_0}(w)=\mbox{proj}_{0}(w)=0 \mbox{ for all } w\in \mathbb{R}^l$, we have by \eqref{formbasic3B} that
\begin{subequations}
 \label{formbasicMC}
 \begin{numcases} { }
   x_i^{k+1/2}=\mbox{proj}_{C_i}(x_i^k-\lambda (B_i^Tv_i^k+  A_i^Tv_{N_1+1}^k)),i=1,2,\cdots,N_1,\label{formbasicMCa}\\
   x_i^{k+1/2}=\mbox{proj}_{C_i}(x_i^k-\gamma\nabla {\theta_i}(x_i^k)-\lambda A_i^Tv_{N_1+1}^k),i=N_1+1,N_1+2,\cdots,N,\label{formbasicMCb}\\
   v_i^{k+1}=(I-\mbox{prox}_{\frac{\gamma}{\lambda}{\theta_i}})(B_ix_i^{k+1/2}+b_i+v_i^k),i=1,2,\cdots,N_1,\label{formbasicMCc}\\
   v_{N_1+1}^{k+1}=\sum_{j=1}^N A_jx_j^{k+1/2}-a+v_{N_1+1}^k,\label{formbasicMCd}\\
   x_i^{k+1}=\mbox{proj}_{C_i}(x_i^k-\lambda (B_i^Tv_i^{k+1}+  A_i^Tv_{N_1+1}^{k+1})),i=1,2,\cdots,N_1,\label{formbasicMCe}\\
   x_i^{k+1}=\mbox{proj}_{C_i}(x_i^k-\gamma\nabla {\theta_i}(x_i^k)-\lambda   A_i^Tv_{N_1+1}^{k+1})),i=N_1+1,N_1+2,\cdots,N,\label{formbasicMCf}
 \end{numcases}
\end{subequations}
where $0<\lambda< 1/(\sum_{i=1}^{N}\lambda_{\max}( A_iA_i^T)+\max \{\lambda_{\max}(B_iB_i^T),i=1,2,\cdots,N_1\})$ and $0<\gamma<2\min\{\beta_i,i=N_1+1,N_1+2,\cdots,N\}$.
It is easy to see that \eqref{formbasicMCa}-\eqref{formbasicMCb}, \eqref{formbasicMCc}-\eqref{formbasicMCd} and \eqref{formbasicMCe}-\eqref{formbasicMCf} can be implemented in  parallel, respectively.
Since \eqref{formbasicMC} can be recast as the original {PDFP} for \eqref{eqbasic_MSB2} which is equal to \eqref{eqbasic_MSB}, we can get the convergence of \eqref{formbasicMC} by the results of  {PDFP}.
Also for some special cases, such as $N_1=0$ and $N_1=N$, one may
even get simpler forms from \eqref{formbasicMC}.
Let $N_1=N$, $B_i=I$ and $b_i=0$ in \eqref{formbasicMC}, we then have
\begin{subequations}
 \label{formbasicMCI}
 \begin{numcases} { }
   x_i^{k+1/2}=\mbox{proj}_{C_i}(x_i^k-\lambda (v_i^k+  A_i^Tv_{N+1}^k)),i=1,2,\cdots,N,\label{formbasicMCIa}\\
   v_i^{k+1}=(I-\mbox{prox}_{\frac{\gamma}{\lambda}{\theta_i}})(x_i^{k+1/2}+b_i+v_i^k),i=1,2,\cdots,N,\label{formbasicMCIc}\\
   v_{N+1}^{k+1}=\sum_{j=1}^{N}A_jx_j^{k+1/2}-a+v_{N+1}^k,\label{formbasicMCId}\\
   x_i^{k+1}=\mbox{proj}_{C_i}(x_i^k-\lambda (v_i^{k+1}+  A_i^Tv_{N+1}^{k+1})),i=1,2,\cdots,N,\label{formbasicMCIe}
  \end{numcases}
\end{subequations}
for solving  \eqref{eqbasic_MSB02}, where $0<\lambda<1/(\sum_{i=1}^{N}\lambda_{\max}( A_iA_i^T)+1)$ and $0<\gamma<+\infty$.
The scheme of \eqref{formbasicMC}, including \eqref{formbasicMCI}, can be implemented in parallel, and there is no requirement for the subproblem solving if the proximity operator of $\theta_i$ have the closed-form representation.

For solving  \eqref{eqbasic_MSB02},  we can also get many others algorithms, by viewing  parts of $\theta_i$ as $f_1$, parts of $\theta_i$ as $f_2\circ B$ and parts of $\theta_i$ as $f_3$. Here we just give an example to show the idea. Let
\begin{align*}
&f_1(x)=0,\  f_2(y)=\chi_0(y),\\
&B=(A_1,A_2,\cdots,A_N),\  b=-a,\  y=Bx+b,\\
&f_3(x)=f_3(x_1,x_2,\cdots,x_N)=\sum_{i=1}^N(\theta_i(x_i)+\chi_{C_i}(x_i)).
\end{align*}
Due to the separability of $f_3$, PDFP \eqref{formbasic3B} can be further expressed as
\begin{subequations}
 \label{formbasicMCI02}
 \begin{numcases} {}
   x_i^{k+1/2}=\underset{x_i \in C_i}{\mbox{argmin  }} \theta_ i(x_i )+\frac {1}{2\gamma} \|x_i -(x_i^k-{\lambda} A_i^T  v^{k})\|^2,i=1,2,\cdots,N,\label{formbasicMCI02a}\\
   v^{k+1}=v^k+(\sum_{j=1}^N A_jx_j^{k+1/2}  -a)\label{formbasicMCI02b},\\
   x_i^{k+1}=\underset{x_i \in C_i}{\mbox{argmin  }} \theta_ i(x_i )+\frac {1}{2\gamma} \|x_i -(x_i^k-{\lambda} A_i^T  v^{k+1})\|^2,i=1,2,\cdots,N,\label{formbasicMCI02c}
 \end{numcases}
\end{subequations}
where $0<\lambda<1/\sum_{i=1}^N\lambda_{\max}( A_iA_i^T)$ and $0<\gamma<+\infty$.
We can write the explicit solution of \eqref {formbasicMCI02a} and \eqref{formbasicMCI02c} for some special $\theta_i$ and $C_i$, for example $\theta_i=\|\cdot\|_1$ and $C_i$ are rectangular domains. If $C_i=\mathbb{R}^{n_i}$, for the schemes \eqref {formbasicMCI02a} and \eqref{formbasicMCI02c}, we just need to work out the proximity operator of $\theta_i$. So the scheme is parallel and easy to implement for solving \eqref{eqbasic_MSB02}, which is the basic problem considered in the context of ADMM.

As shown in Section 1, we can write the problem \eqref{eqbasic_MSB} (or problem \eqref{eqbasic_TVL1}) in the form \eqref{PDFP2O3B:eqbasic} with many other ways, and then derive new schemes to solve it in terms of PDFP \eqref{formbasic3B}. Since the discussion is routine, we omit the details. What we have to emphasize is that our method for constructing algorithms for solving problem \eqref{eqbasic_MSB} or \eqref{eqbasic_TVL1} is very flexible.

\subsection{Comparison to ADMM-like algorithms}
In this subsection, let us show the difference of ADMM and PDFP for \eqref{eqbasic_MSB} by solving the following problem:
\begin{equation}
\label{ThreeB:eqbasic_LC}
\begin{aligned}
   &{\mbox{min}}\quad \theta_ 1(x_1 )+\theta_ 2(x_2 ) +\theta_ 3(x_3 )\\
   &\mbox{st. }A_1x_1 +A_2x_2  +A_3x_3=a,\\
   &\phantom {\mbox{st. }}x_1 \in C_1, x_2 \in C_2, x_3 \in C_3,
\end{aligned}
\end{equation}
where $x_i \in \mathbb{R}^{n_i},i=1,2,3$.

For the ADMM method, the above problem is first transformed to solve the following min-max problem:
\begin{align}
   \underset{x_1 \in C_1, x_2 \in C_2,x_3\in C_3}{\mbox{ min}}\underset{w}{\mbox{ max}}\quad \mathcal{L}_{\beta}(x_1, x_2, x_3, w)=\sum_{i=1}^3\theta_ i(x_i ) -\langle w,\sum_{i=1}^3 A_ix_i -a\rangle+\frac \beta 2\|\sum_{i=1}^3 A_ix_i -a\|^2. \label{ThreeB:eqbasic_AL}
\end{align}
Let $v=w/\beta$. We then use the alternating direction method to solve problem \eqref{ThreeB:eqbasic_AL}, leading to  the following algorithm
\begin{subequations}
 \label{ThreeB_ADMM}
 \begin{numcases} { }
   x_1 ^{k+1}=\underset{x_1 \in C_1}{\mbox{argmin  }} \theta_ 1(x_1 )+\frac {\beta} 2\|A_1x_1 +(A_2x_2^k+A_3x_3^k-v^k-a)\|^2,\label{ThreeB_ADMMa}\\
   x_2^ {k+1}=\underset{x_2 \in C_2}{\mbox{argmin  }} \theta_ 2(x_2 )+\frac {\beta} 2\|A_2x_2 +(A_1x_1^{k+1}+A_3x_3^k-v^k-a)\|^2,\label{ThreeB_ADMMb}\\
   x_3^ {k+1}=\underset{x_3 \in C_3}{\mbox{argmin  }} \theta_ 3(x_3 )+\frac {\beta} 2\|A_3x_3 +(A_1x_1^{k+1}+A_2x_2^{k+1}-v^k-a)\|^2,\label{ThreeB_ADMMc}\\
   v^{k+1}=v^k- \tau (A_1x_1 ^{k+1}+A_2x_2^ {k+1}+A_3x_3^{k+1}-a),\label{ThreeB_ADMMd}
  \end{numcases}
\end{subequations}
In general, \eqref{ThreeB_ADMMa}-\eqref{ThreeB_ADMMc} need to solve three subprograms whenever $A_i\neq I$ and the scheme is not a parallel algorithm. In addition, if one of \eqref{ThreeB_ADMMa}-\eqref{ThreeB_ADMMb} is not easy to solve due to the constraints $C_i$, we must introduce new auxiliary variables to get the solution. Though the treatment is routine, the solution process will become rather complicated. More importantly, as showed in \cite{CHY14}, the scheme \eqref{ThreeB_ADMM} is not necessarily convergent if there is no further assumption on \eqref{ThreeB:eqbasic_LC}.  Recently it is popular to propose some variants of ADMM to overcome this disadvantage, for example, some prediction-correction methods were proposed in \cite{HTY12}, and the Jacobian decomposition of augmented Lagrangian method (ALM) with proximal terms was introduced in \cite{HXY13}.

Now, let us continue to show how to solve \eqref{ThreeB:eqbasic_LC} in view of PDFP. By using indicator functions, \eqref{ThreeB:eqbasic_LC} is equivalent to
\begin{equation}
   {\mbox{min }} \sum_{i=1}^3\theta_i(x_i)+\chi_{0}(\sum_{i=1}^3 A_ix_i-a) +\sum_{i=1}^3\chi_{C_i}(x_i). \label{ThreeB:eqbasic_PDFP}
\end{equation}
Then we can use PDFP to solve \eqref{ThreeB:eqbasic_PDFP} in various forms. For example, by setting $N=3$ in \eqref{formbasicMCI02}, we can get the following algorithm
\begin{subequations}
 \label{ThreeB:form3B04}
 \begin{numcases} {}
   x_i^{k+1/2}=\underset{x_i \in C_i}{\mbox{argmin  }} \theta_ i(x_i )+\frac {1}{2\gamma} \|x_i -(x_i^k-{\lambda} A_i^T  v^{k})\|^2,i=1,2,3,\label{ThreeB:form3B04a}\\
   v^{k+1}=v^k+(A_1x_1^{k+1/2} +A_2x_2^{k+1/2}+A_3x_3^{k+1/2} -a)\label{ThreeB:form3B04b},\\
   x_i^{k+1}=\underset{x_i \in C_i}{\mbox{argmin  }} \theta_ i(x_i )+\frac {1}{2\gamma} \|x_i -(x_i^k-{\lambda} A_i^T  v^{k+1})\|^2,i=1,2,3,\label{ThreeB:form3B04c}
 \end{numcases}
\end{subequations}
where $0<\lambda<1/\sum_{i=1}^{3}\lambda_{\max}( A_iA_i^T)$ and $0<\gamma<+\infty$. Compared to the scheme of \eqref{ThreeB_ADMM}, the scheme of \eqref{ThreeB:form3B04} is parallel and always convergent. Nevertheless, the computation cost increases with the addition of a symmetric step, which may double the work of each step.
To avoid the disadvantage, we can also extend the scheme in \cite{CP11,C13,LSXZ15,LZ15} with the same treatment given above.

When the subproblems in \eqref{ThreeB:form3B04a} are not easy to solve due to the constraints $C_i$, we can also use  \eqref{formbasicMCI} and  get
\begin{subequations}
 \label{ThreeB:form3B03}
 \begin{numcases} {}
   x_i^{k+1/2}=\mbox{proj}_{C_i}(x_i^k-{\lambda} (v_i^k+A_i^T v_4^{k})),i=1,2,3,\label{ThreeB:form3B03a}\\
   v_i^{k+1}=(I-\mbox{prox}_{\frac{\gamma}{\lambda}{{\theta_i}}})(x_i^{k+1/2}+v_i^k),i=1,2,3,\label{ThreeB:form3B03b}\\
   v_4^{k+1}=v_4^k+(A_1x_1^{k+1/2} +A_2x_2^{k+1/2}+A_3x_3^{k+1/2} -a),\label{ThreeB:form3B03b03}\\
   x_i^{k+1}=\mbox{proj}_{C_i}(x_i^k-{\lambda} (v_i^k+A_i^T v_4^{k+1})),i=1,2,3,\label{ThreeB:form3B03c}
\end{numcases}
\end{subequations}
where $0<\lambda<1/(\sum_{i=1}^{3}\lambda_{\max}( A_iA_i^T)+1)$ and $0<\gamma<+\infty$.

If $\theta_i$ are both differentiable  with  $1/\beta_i$-Lipschitz continuous gradient, respectively.
We can set $N_1=0$ and $N=3$ in \eqref{formbasicMC} to get an furtherly linearized scheme as
\begin{subequations}
 \label{ThreeB:form3B01}
 \begin{numcases} {}
   x_i^{k+1/2}=\mbox{proj}_{C_i}(x_i^k-\gamma\nabla {{\theta_i}}(x_i^k)-{\lambda} A_i^T  v^{k}),i=1,2,3,\label{ThreeB:form3B01a}\\
   v^{k+1}=v^k+(A_1x_1^{k+1/2} +A_2x_2^{k+1/2}+A_3x_3^{k+1/2} -a)\label{ThreeB:form3B01b},\\
   x_i^{k+1}=\mbox{proj}_{C_i}(x_i^k-\gamma\nabla {{\theta_i}}(x_i^k)-{\lambda} A_i^T  v^{k+1}),i=1,2,3,\label{ThreeB:form3B01c}
 \end{numcases}
\end{subequations}
where  $0<\lambda<1/\sum_{i=1}^{3}\lambda_{\max}( A_iA_i^T)$ and $0<\gamma<2\min\{\beta_1,\beta_2,\beta_3\}$.

\section{Numerical experiments} \label{sec:numericalE}

In this section, we will illustrate the application of  {PDFP} for multi-block problems  through two examples, related to \eqref{eqbasic_TVL1} and \eqref{eqbasic_MSB}, respectively. The first one is the total variation regularized computerized tomography (CT) reconstruction with constraints, and the second one is on some quadratic programming or linear equation examples given in \cite{CHY14} as the counter examples for the convergence of thee-block ADMM.

\subsection{CT reconstruction}

The standard CT reconstruction algorithm in clinical applications is the so-called Filtered Back Projection (FBP) algorithm.
In the presence of noise, this problem becomes difficult since the inverse of Radon transform is unbounded and ill-posed. In the literature, the model is constructed based on TV regularization \eqref{eq_TVL2}, i.e
\begin{align*}
   x^*=\underset{x\in {C}}{\mbox{arg min}}\quad \frac 1 2\|Ax-a\|^2+\mu \|Dx\|_1.
\end{align*}
Here $A$ is the Radon transform matrix, $a$ is the measured projections vector, and $D$ is the discrete gradient operator.  The size of $A$ is generally huge and it is very difficult for us to efficiently solve a linear system with $A$ as the coefficient matrix. $\|Dx\|_{1}$ is  the usual $\ell_1$ based regularization in order to promote sparsity under the transform $D$ and $\mu>0$ is the regularization parameter. To be more precise, we use the isotropic total variation as the regularization term,  and assume that the solution should belong to [0,255], in other words, the constraint set is defined as  $C=\{x=(x_1,x_2,\cdots, x_n)^T\in \mathbb{R}^n|x_i\in [0,255],i=1,2,\cdots,n\}$.  We have shown in \cite{CHZ1302} that it is useful to impose the above constraints in CT to improve the quality of reconstructed images.

In our numerical simulation, we still use the same example tested in \cite{ZBO11}, i.e.,  $50$ uniformly oriented projections are simulated for a
$128  \times 128$ Shepp-Logan phantom image and then white Gaussian noise of  mean $0$ and variance $1$ is added to the data.  For this example,
we compute $\lambda_{\max}(AA^T)=1.5086$. It is well known in total variation application that $\lambda_{\max}(DD^T)=8$. So we set $\gamma=1.3$, $\lambda=1/8$ ($0<\gamma<2/\lambda_{\max}(AA^T)=1.3257$ and $0<\lambda<1/8$ in PDFP according to Theorem 3.1 in \cite{CHZ15} in Scheme 1 (cf. \eqref{formTVL2C2}). Correspondingly we set $\gamma=1.3$, $\lambda=1/9$ in Scheme 2 (cf. \eqref{formTVL2C1}).  Set $\gamma=20$ and $\lambda=1/(8+1.5086)$ in Scheme 3  (cf. \eqref{formTVL2nof102}). Set $\gamma=100$, $\lambda=1/(8+N)$, $N=20$ in Scheme 4 (cf. \eqref{formTVL2nof103} and \eqref{keytechnique02}). Here we do not implement \eqref{formTVL2nof1} since it needs to solve a large linear system, nor \eqref{formTVL2nof104} as it is a parallel form of Scheme 3.

\begin{figure}[!htp]\centering
\caption{$log_{10}$(energy) and  PSNR \textit{versus} iterations for different PDFP algorithms in CT reconstruction.}
 \label{figure:CT_cmp}
 \begin{tabular}{cc}
    \includegraphics[width=0.4 \textwidth]{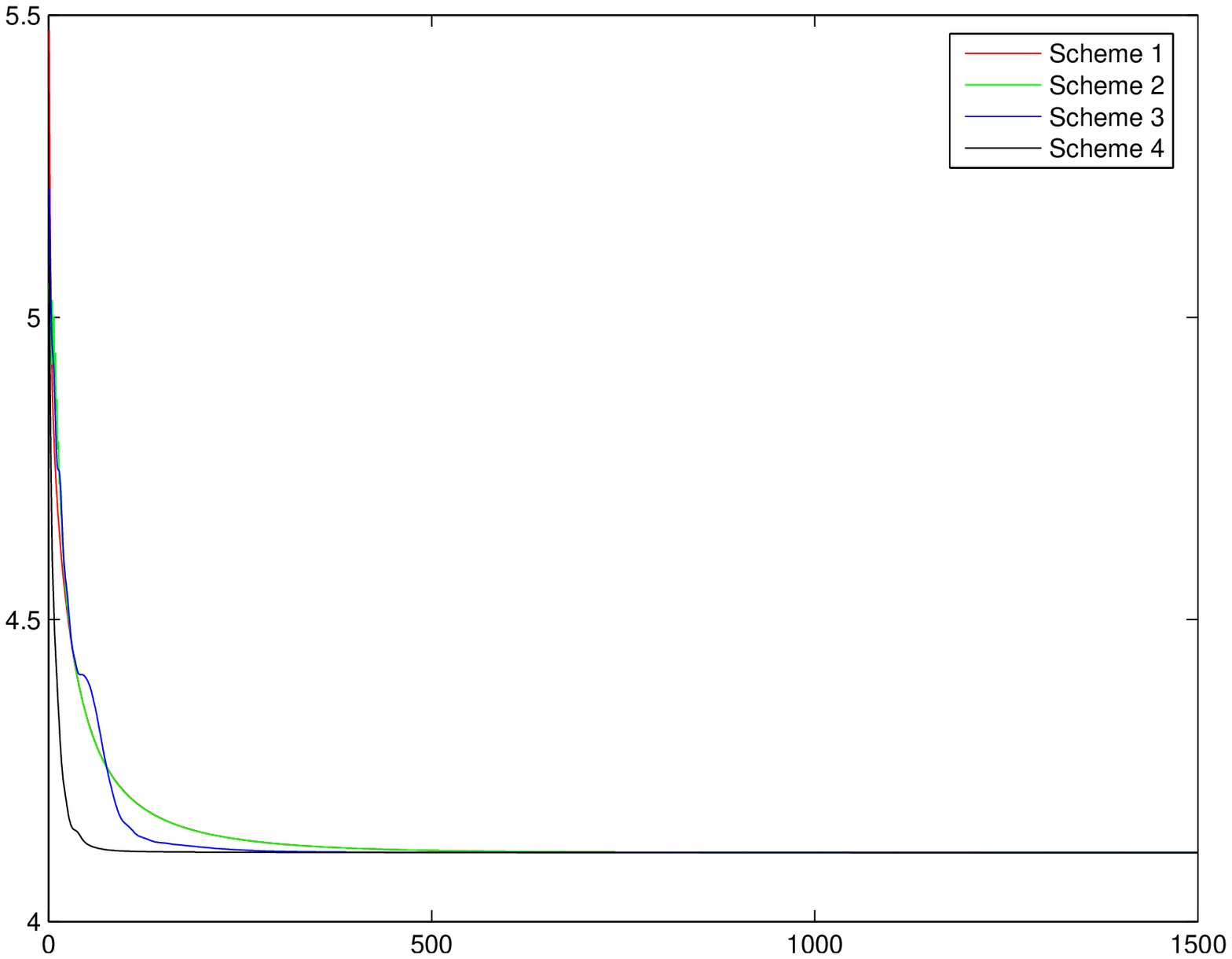}&
    \includegraphics[width=0.4 \textwidth]{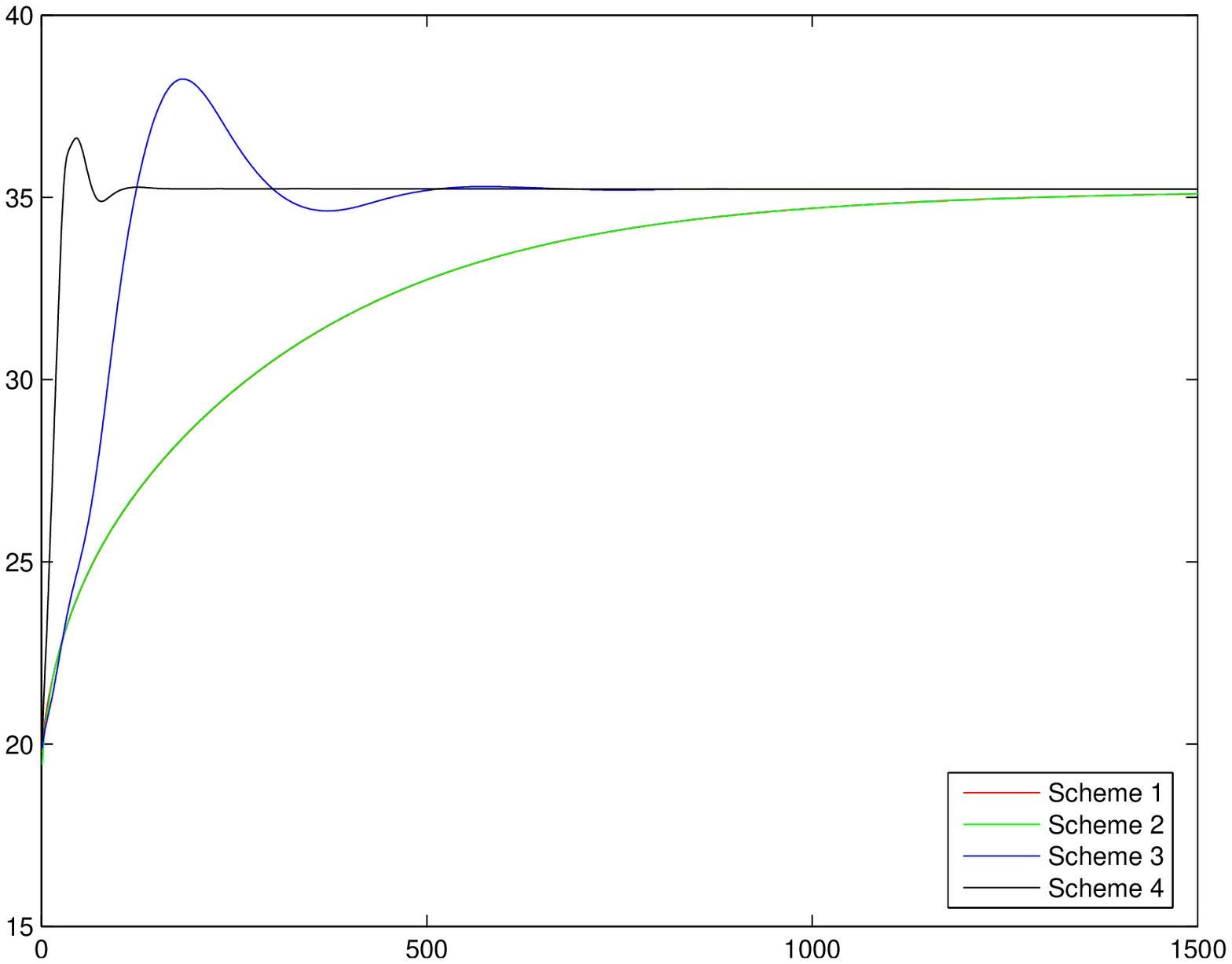}
 \end{tabular}
\end{figure}

\begin{figure}[!htp]\centering\small
\caption{The best recovery results for CT in 1500 iterations. }
 \label{figure:CT_best_cmp}
 \begin{tabular}{cccccc}
   &\small Scheme 1 &\small Scheme 2 &\small Scheme  3& \small Scheme 4   \\
   &\includegraphics[width=0.2\textwidth]{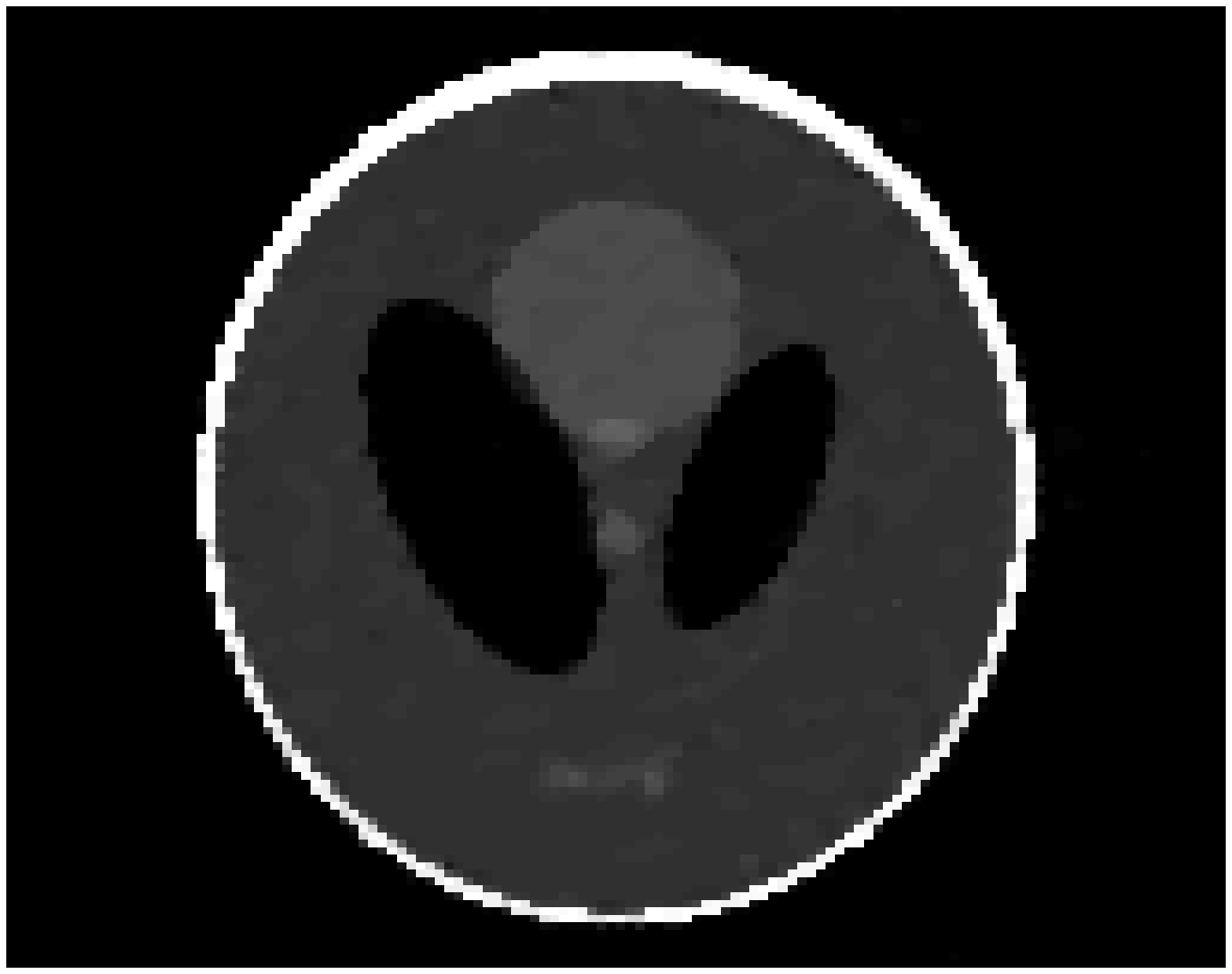}&
   \includegraphics[width=0.2\textwidth]{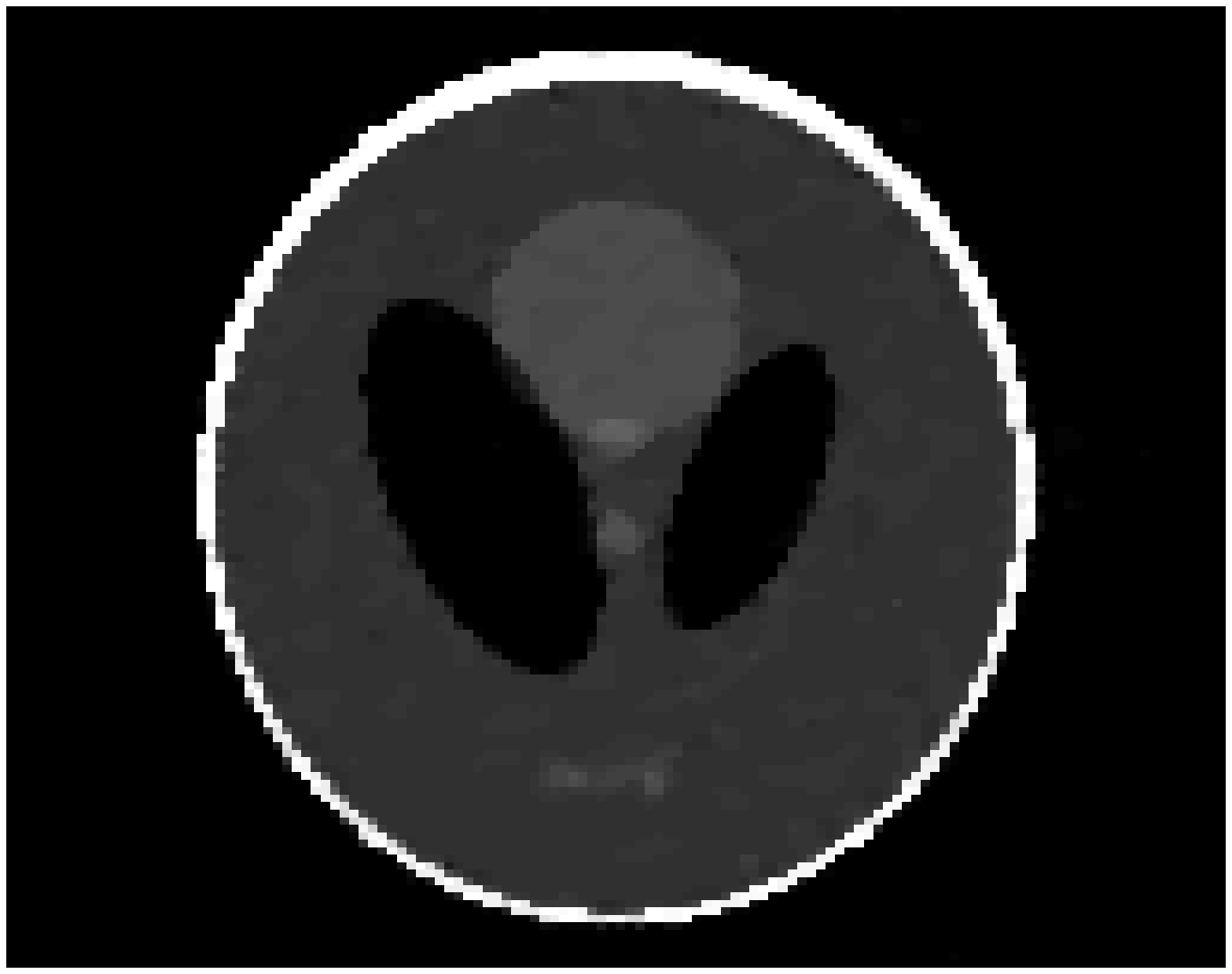}&
   \includegraphics[width=0.2\textwidth]{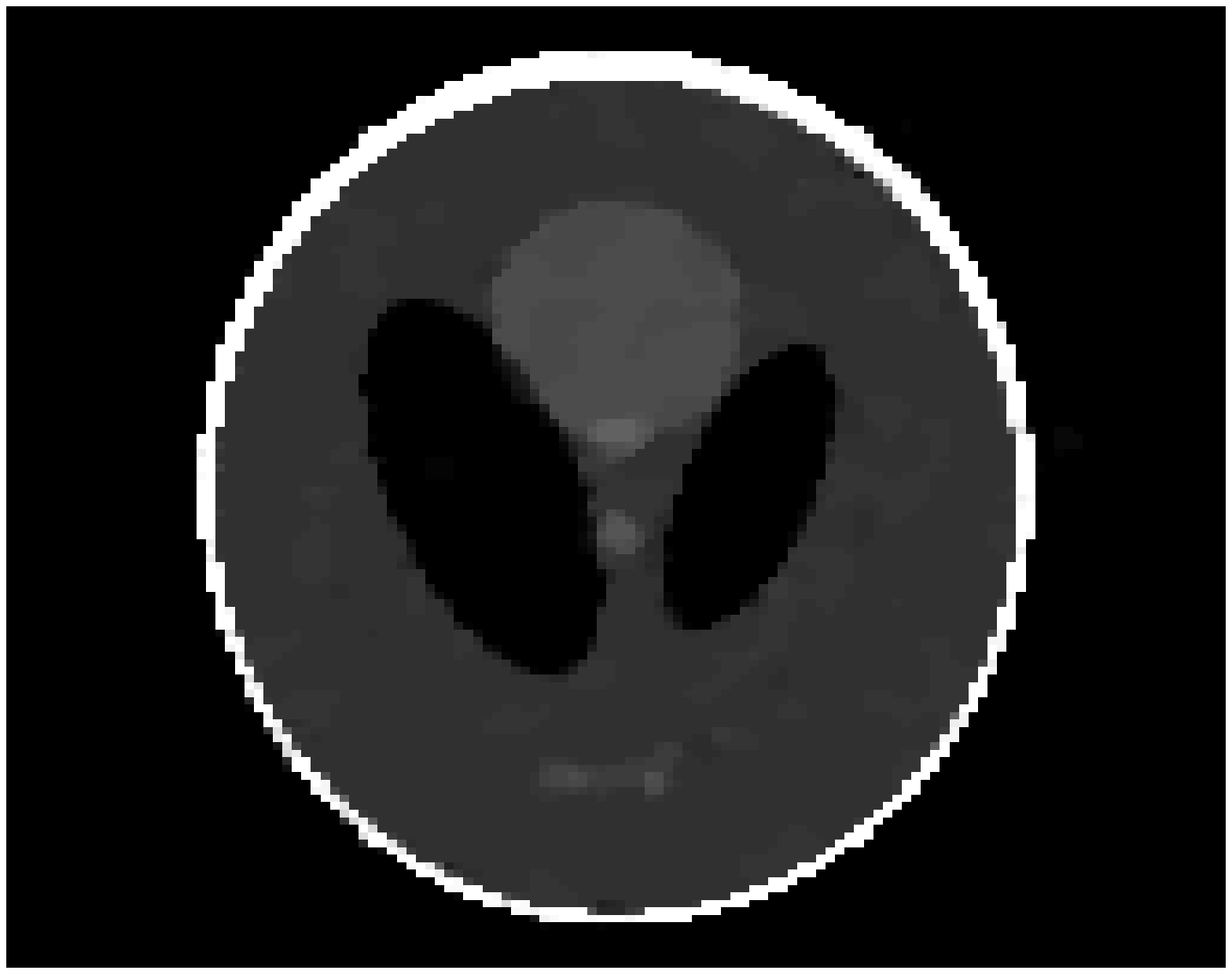}&
   \includegraphics[width=0.2\textwidth]{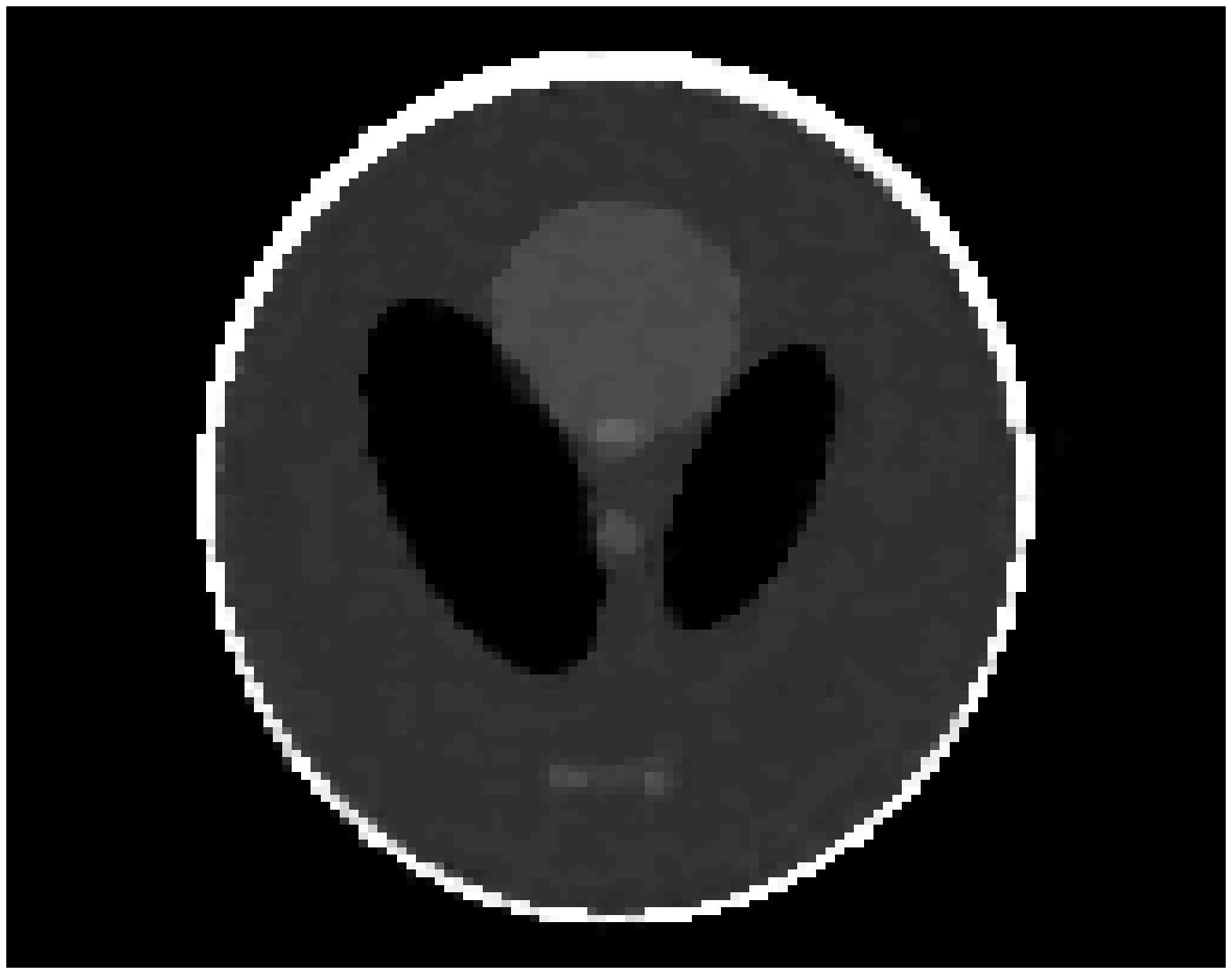}
   \\
   itn&1500&1500&183&45\\
   time&35.33&   34.76 &   4.29  & 23.53\\
   PSNR&35.0997&   35.1003 &  38.2467  & 36.6266
 \end{tabular}
\end{figure}

From Figure \ref{figure:CT_cmp}, we can see that Scheme 3 and Scheme 4 can get relatively better results with higher PSNR and  use far less iteration steps than Scheme 1 and Scheme 2. According to Theorem 3.2 in \cite{CHZ15}, the convergence rate of PDFP depends on the Lipschitz constant of $\nabla f_1$ (the smaller the better)
and the quantity $\delta$ indicating the strongly monotone nature of $\partial f_2^*$ (the larger the better).
So an intuitive explanation for our previous observation is that the related Lipschitz constant of the
gradient of the  function $\frac 1 2\|Ax-a\|^2$, the largest eigenvalue of $A^TA$, is relatively large, which implies the slow
convergence of Scheme 1 and Scheme 2. On the other hand, if we view $\frac 1 2\|Ax-a\|^2$ as a part of  $f_2\circ B$, then $f_1$ is taken to be $zero$ and the corresponding parameter $\delta$ of $\partial f_2^*$ become larger, which would thus improve the convergence rate of the algorithm. The problems in Scheme 3  and Scheme 4 are how to choose the arbitrary parameter $\gamma$ so that we can get faster convergence.
In addition, Scheme 3  and Scheme 4  have a relatively high PSNR in the first steps and then keep the results almost unchanged. The best PSNR of Scheme 3 are better than the ones in Scheme 4. Scheme 3 requires a little more steps than Scheme 4 but the computation time is far less. The times shown in Figure \ref{figure:CT_best_cmp} are the ones when the underlying algorithms are carried out in a sequential way.
As a matter of fact, all the schemes proposed here can be implemented in parallel with ease, which will reduce the computation time essentially.

\subsection{Application to non convergent examples for the direct extension of ADMM}\label{ADMM_notConvergence}
\begin{figure}[!htp]\centering
\caption{Errors  \textit{versus} iterations for  PDFP within 2000 steps. }
 \label{figure:ADMMnotConvergence}
 \begin{tabular}{ccc}
    \includegraphics[width=0.3\textwidth]{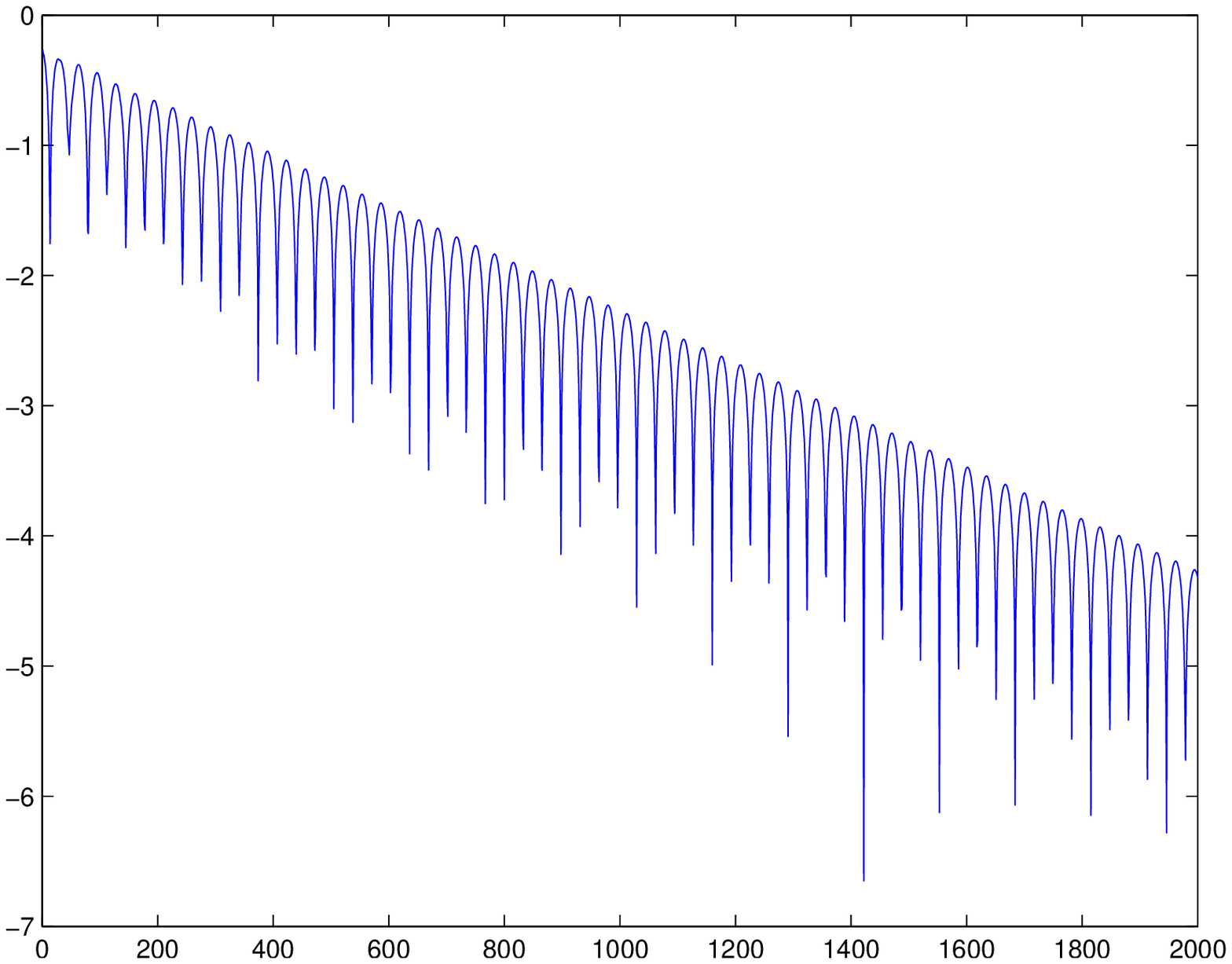}&
    \includegraphics[width=0.3\textwidth]{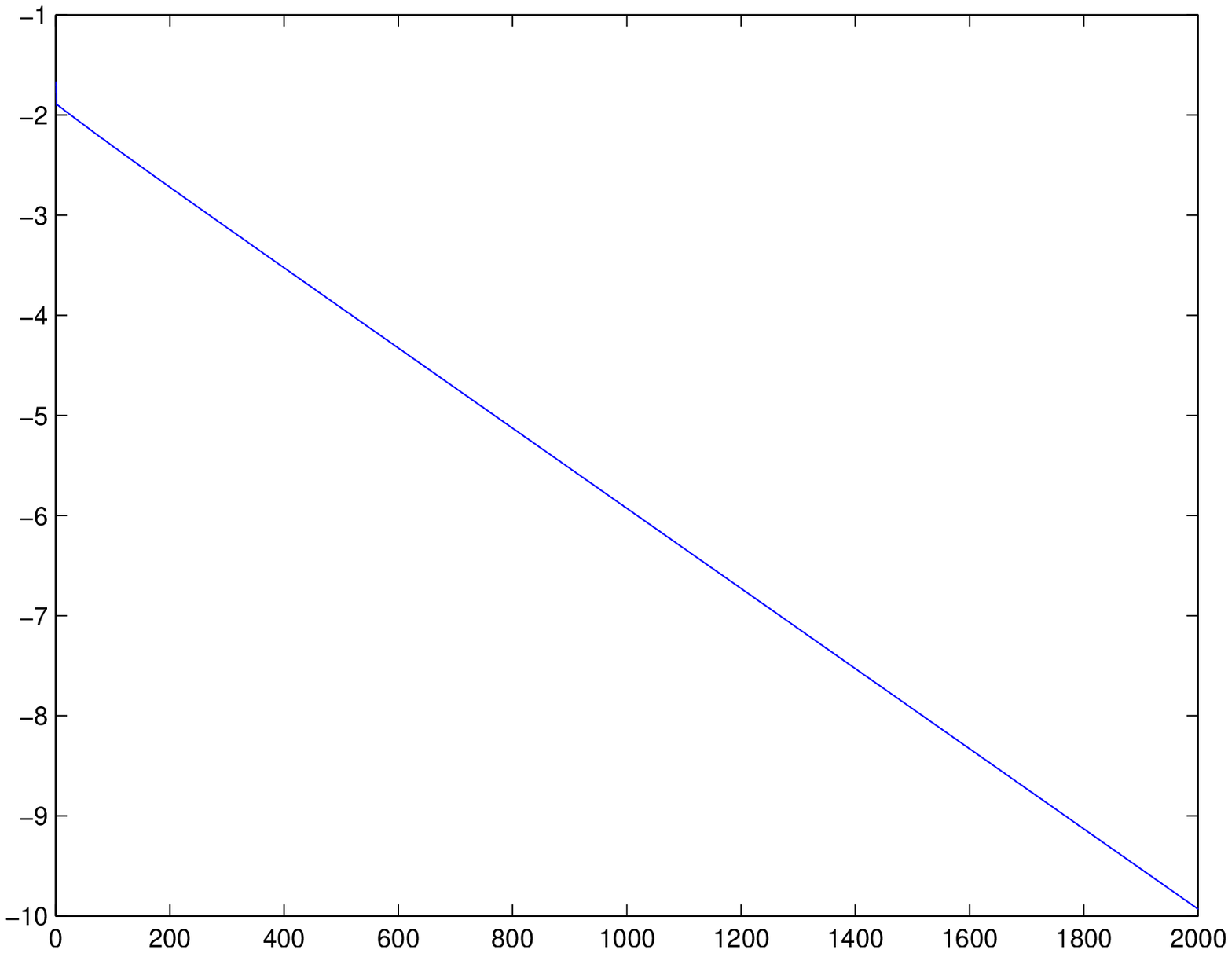}&
    \includegraphics[width=0.3\textwidth]{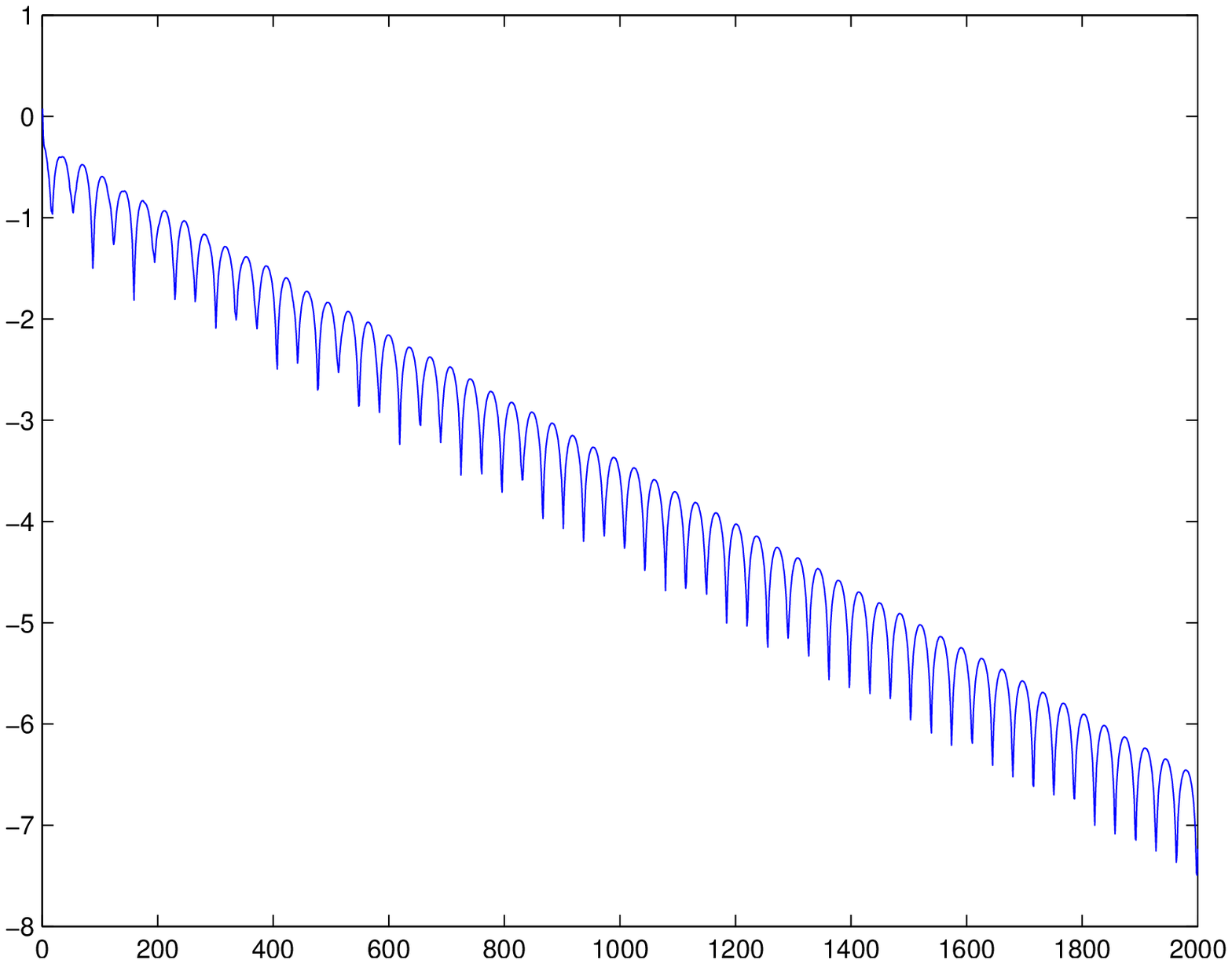}\\
    (a) Example 1&(b) Example 2&(c) Example 3
 \end{tabular}
\end{figure}

As showed in \cite{CHY14}, the direct extension scheme \eqref{ThreeB_ADMM} is not necessarily convergent if there is no further assumption on \eqref{ThreeB:eqbasic_LC}. Some non-convergent examples of ADMM are given in \cite{CHY14}. We will use these simple but important examples to test the properties of our PDFP schemes. We know that they are convergent by the theory of PDFP developed in \cite{CHZ15}.
Thus this provides alternative approach when ADMM does not converge for some applications. The errors with respect to the true solution within 2000 steps are given in Figure \ref{figure:ADMMnotConvergence}.

The first example is solving linear equation
\begin{align}
   \begin{pmatrix}
   1\\
   1\\
   1
   \end{pmatrix}x_1 +\begin{pmatrix}
   1\\
   1\\
   2
   \end{pmatrix}x_2 +\begin{pmatrix}
   1\\
   2\\
   2
   \end{pmatrix}x_3=\begin{pmatrix}
   0\\
   0\\
   0
   \end{pmatrix}.\label{Ax0:eqbasic_LC}
\end{align}
\eqref{Ax0:eqbasic_LC} is a special case of \eqref{ThreeB:eqbasic_LC}, where $A=(A_1,A_2,A_3)=\begin{pmatrix}
   1& 1& 1\\
   1 &1 &2\\
   1 &2 &2
   \end{pmatrix}
   $ and
$a=\begin{pmatrix}
   0\\
   0\\
   0
   \end{pmatrix}$.
It is easy to verify that $A$ is nonsingular, and the true solution is $x_1=0$, $x_2=0$ and $x_3=0$. Moreover, the corresponding optimal Lagrange multipliers are all 0.  Let $\theta_i=0$ and $C_i=\mathbb{R}$, i=1,2,3 in \eqref{ThreeB:form3B01}, \eqref{ThreeB:form3B03} or \eqref{ThreeB:form3B04}, we can get the following scheme to solve it.
Namely
\begin{subequations}
 \label{Ax:form3B01}
 \begin{numcases} {}
   x_i^{k+1/2}=x_i^k-{\lambda} A_i^T  v^{k},i=1,2,3,\label{Ax:form3B01a}\\
   v^{k+1}=v^k+(A_1x_1^{k+1/2} +A_2x_2^{k+1/2}+A_3x_3^{k+1/2} -a),\label{Ax:form3B01b}\\
   x_i^{k+1}=x_i^k-{\lambda} A_i^T  v^{k+1},i=1,2,3,\label{Ax:form3B01c}
 \end{numcases}
\end{subequations}
where $0<\lambda\leq 1/\sum_{i=1}^{3}\lambda_{\max}( A_i^TA_i)$. Substitute $x_i^{k+1/2}$ with $x_i^{k+1}$, and \eqref{Ax:form3B01} implies
\begin{subequations}
 \label{Ax:form3B02}
 \begin{numcases} {}
   x_i^{k+1}=x_i^k-\lambda A_i^T(A_1x_1^{k} +A_2x_2^{k}+A_3x_3^{k} -a)-{\lambda} A_i^T  v^{k},i=1,2,3,\label{Ax:form3B02a}\\
   v^{k+1}=v^k+(A_1x_1^{k+1} +A_2x_2^{k+1}+A_3x_3^{k+1} -a),\label{Ax:form3B02b}
 \end{numcases}
\end{subequations}
i.e.,
\begin{align}
   x_i^{k+1}=x_i^k-\lambda A_i^T \sum_{j=1}^k(A_1x_1^{j} +A_2x_2^{j}+A_3x_3^{j} -a)-{\lambda} A_i^T  v^{0},i=1,2,3,
\end{align}
We set $\lambda= 1/\sum_{i=1}^{3}\lambda_{\max}( A_i^TA_i)$ and the initial values of the elements of $x_i$ and $v$ as $1$.

The second example is solving
\begin{align}
    &{\mbox{min}}\quad 0.05x_1^2+0.05x_2^2 +0.05x_3^2\nonumber\\
   &\mbox{st. }\begin{pmatrix}
   1\\
   1\\
   1
   \end{pmatrix}x_1 +\begin{pmatrix}
   1\\
   1\\
   2
   \end{pmatrix}x_2 +\begin{pmatrix}
   1\\
   2\\
   2
   \end{pmatrix}x_3=\begin{pmatrix}
   0\\
   0\\
   0
   \end{pmatrix}.\label{x2Ax:eqbasic_LC}
\end{align}
\eqref{x2Ax:eqbasic_LC} is also a special case of \eqref{ThreeB:eqbasic_LC}.
Let $\theta_i=0.05x_i^2$ and $C_i=\mathbb{R}$, $i=1,2,3$ in \eqref{ThreeB:form3B01}, we can easily get the algorithm for solving \eqref{x2Ax:eqbasic_LC}, where $0<\lambda\leq 1/\sum_{i=1}^{3}\lambda_{\max}( A_i^TA_i)$ and $0<\gamma<2/0.1=20$. According to the  convergence rate theory about PDFP$^2$O given in \cite{CHZ13},  this algorithm has linear convergence rate, which is also confirmed by Figure  \ref{figure:ADMMnotConvergence}(b) , since $f_1(x)=f_1(x_1,x_2,x_3)=\sum_{i=1}^30.05x_i^2$ is strongly convex and $BB^T=AA^T$ is positive symmetric definite in \eqref{PDFP2O3B:eqbasic} with $f_3=0$. Set $\gamma=1/0.1=10$ and the others setting are the same as the first example.

The third example given in \cite{CHY14} is more sophisticated. It can be described as
\begin{align}
    &{\mbox{min}}\quad 0.5x_1^2\nonumber\\
   &\mbox{st. }\begin{pmatrix}
   1\\
   1\\
   1
   \end{pmatrix}x_1+\begin{pmatrix}
   1\\
   1\\
   1
   \end{pmatrix}x_2 +\begin{pmatrix}
   1\\
   1\\
   2
   \end{pmatrix}x_3 +\begin{pmatrix}
   1\\
   2\\
   2
   \end{pmatrix}x_4=\begin{pmatrix}
   0\\
   0\\
   0
   \end{pmatrix}.\label{x12Ax:eqbasic_LC}
\end{align}
The feasible region of \eqref{x12Ax:eqbasic_LC} is not a singleton, and the objective function is only related with $x_1$. The optimal solution of \eqref{x12Ax:eqbasic_LC} is $x_i=0$, $i=1,2,3,4$. Similar to \eqref{ThreeB:form3B01}, let $N_1=0$, $N=4$, $\theta_1=0.5x_1^2$, $\theta_i=0$, $i=2,3,4$, $C_i=\mathbb{R}$, $i=1,2,3,4$ in \eqref{formbasicMC}, we can easy to get the algorithm for solving \eqref{x12Ax:eqbasic_LC}, where $0<\lambda\leq1/\sum_{i=1}^{4}\lambda_{\max}( A_i^TA_i)$ and $0<\gamma<2$. We set $\lambda=1/\sum_{i=1}^{4}\lambda_{\max}( A_i^TA_i)$ and $\gamma=1$. The others setting are same as the first example. From Figure  \ref{figure:ADMMnotConvergence}, we can see that PDFP solves these problems with linear convergence rate, which are divergent examples using the direct extension ADMM \eqref{ThreeB_ADMM}.

\section{Conclusion}
We extend the ideas of a proximal primal-dual fixed point algorithm {PDFP} to solve separable multi-block minimization problems with and without linear constraints.
The variants of {PDFP} are fully decoupled and therefore easy to implement. Moreover, the algorithms are parallel naturally, so they are very suitable for solving   large-scale problems from real-world models.
Through numerical experiments, we can see that treating smooth functions as parts of $f_2\circ B$ leads to better convergence and partial inverse can be viewed as preconditioner for a good balance of convergence speed and computational cost. The convergence conditions on the parameters are arbitrary positive numbers, while the choices might heavily affect the convergent speed, which may make parameter choosing a difficult problem in practice. Therefore the proper decomposition of the smooth functions and non-smooth functions, and explicit or implicit schemes should depend on the  properties and computation balances in real applications. Finally, for problems with constraints that three block ADMM may fail to converge,
PDFP algorithm can be also a choice with the guarantee of theoretical convergence.

\medskip

\noindent{\bf Acknowledgements.} P. Chen was partially supported by the PhD research startup foundation of Taiyuan University of Science and Technology (No. 20132024). J. Huang was partially supported by NSFC (No. 11571237). X. Zhang was  partially supported by NSFC (No. 91330102 and GZ1025) and 973 program (No. 2015CB856004).

\end{document}